\documentclass[11pt]{article}

\usepackage{amscd,amsmath,amssymb}

\def\eqref#1{(\ref{#1})}
\newcommand{\goth}{\frak}

\newcommand{\arrow}{{\:\longrightarrow\:}}

\newcommand{\C}{{\Bbb C}}
\newcommand{\R}{{\Bbb R}}

\newcommand{\1}{\sqrt{-1}\:}

\renewcommand{\c}[1]{{\cal #1}}
\newcommand{\calo}{{\cal O}}

\renewcommand{\bar}{\overline}
\renewcommand{\phi}{\varphi}
\renewcommand{\epsilon}{\varepsilon}
\renewcommand{\geq}{\geqslant}
\renewcommand{\leq}{\leqslant}

\newcommand{\Vol}{\operatorname{Vol}}

\newcommand{\comment}[1]{{}}

\def\blacksquare{\hbox{\vrule width 4pt height 4pt depth 0pt}}
\def\endproof{\blacksquare}

\makeatletter

\@ifundefined{Bbb}
     {\newcommand{\Bbb}[1]{{\mathbb #1}}}%
{}

\newcommand{\ps@verbit}{%
  \renewcommand{\@oddhead}{%
          \scriptsize
          {Wirtinger numbers}
          \hfil\tiny {M. Verbitsky, \ \ \ December 7, 1998}}
  \renewcommand{\@evenhead}{\@oddhead}
  \renewcommand{\@oddfoot}{\hfil\thepage\hfil}
  \renewcommand{\@evenfoot}{\@oddfoot}}
 
\newcounter{Mycounter}[section]
\newcounter{lemma}[section]
\renewcommand{\thelemma}{{Lemma \thesection.\arabic{lemma}}}
\newcommand{\lemma}{%
     \setcounter{lemma}{\value{Mycounter}}
     \refstepcounter{lemma}
     \stepcounter{Mycounter}
     {\bf \thelemma:\ }}

\newcounter{claim}[section]
\newcounter{sublemma}[section]
\newcounter{corollary}[section]
\renewcommand{\thecorollary}{{Corollary \thesection.\arabic{corollary}}}
\newcommand{\corollary}{%
     \setcounter{corollary}{\value{Mycounter}}
     \refstepcounter{corollary}
     \stepcounter{Mycounter}
     {\bf \thecorollary:\ }}

\newcounter{theorem}[section]
\renewcommand{\thetheorem}{{Theorem \thesection.\arabic{theorem}}}
\newcommand{\theorem}{%
     \setcounter{theorem}{\value{Mycounter}}
     \refstepcounter{theorem}
     \stepcounter{Mycounter}
     {\bf \thetheorem:\ }}

\newcounter{conjecture}[section]
\newcounter{proposition}[section]
\renewcommand{\theproposition}
       {{Proposition \thesection.\arabic{proposition}}}
\newcommand{\proposition}{%
     \setcounter{proposition}{\value{Mycounter}}
     \refstepcounter{proposition}
     \stepcounter{Mycounter}
     {\bf \theproposition:\ }}

\newcounter{definition}[section]
\renewcommand{\thedefinition}
       {{Definition \thesection.\arabic{definition}}}
\newcommand{\definition}{%
     \setcounter{definition}{\value{Mycounter}}
     \refstepcounter{definition}
     \stepcounter{Mycounter}
     {\bf \thedefinition:\ }}

\newcounter{example}[section]
\newcounter{remark}[section]
\renewcommand{\theremark}{{Remark \thesection.\arabic{remark}}}
\newcommand{\remark}{%
     \setcounter{remark}{\value{Mycounter}}
     \refstepcounter{remark}
     \stepcounter{Mycounter}
     {\bf \theremark:\ }}

\newcounter{problem}[section]
\newcounter{question}[section]
\renewcommand{\thequestion}{{Question \thesection.\arabic{question}}}
\newcommand{\question}{%
     \setcounter{question}{\value{Mycounter}}
     \refstepcounter{question}
     \stepcounter{Mycounter}
     {\bf \thequestion:\ }}

\@addtoreset{equation}{section}
\@addtoreset{footnote}{section}
\makeatother

\begin{document}

\begin{center}
{\Large\bf
Wirtinger numbers\\ and
holomorphic symplectic immersions}
\\[4mm]
Misha Verbitsky,\\[4mm]
{\tt verbit@@thelema.dnttm.rssi.ru}
\end{center}

{\small 
\hspace{0.2\linewidth}
\begin{minipage}[t]{0.7\linewidth}
For any subvariety of a compact holomorphic symplectic K\"ahler 
manifold, we define the symplectic Wirtinger number $W(X)$. We show that
$W(X)\leq 1$, and the equality is reached if and only if 
the subvariety $X\subset M$ is trianalytic, i. e. compactible
with the hyperk\"ahler structure on $M$. For a sequence 
$X_1 \arrow X_2 \arrow ... X_n\arrow M$
of immersions of simple holomorphic symplectic manifolds,
we show that $W(X_1) \leq W(X_2) \leq ... \leq W(X_n)$.
\end{minipage}
}

\tableofcontents


\section{Introduction}


It is well known that every
hyperk\"ahler manifold is holomorphic symplectic
(see \eqref{_holomo_symple_Equation_}).
In other words, we have a natural functor $\c F$
from the category of hyperk\"ahler manifolds (morphisms -- all immersions
compatible with metric and quaternionic action) to the category
of holomorphic symplectic manifolds (morphisms -- all symplectic immersions).
As follows from Yau's solution of the 
Calabi conjecture (\ref{_Calabi_Yau_Theorem_}),
every compact holomorphic symplectic K\"ahler 
manifold admits a hyperk\"ahler metric, which is uniquely
determined by its K\"ahler class. Thus, $\c F$
establishes a bijection on objects of the appropriate
categories. However, $\c F$ is not an equivalence of categories,
in any reasonable sense: there are holomorphic symplectic
immersions which are not compatible with a hyperk\"ahler structure.
One of such examples was found in \cite{_Verbitsky:Hilbert_}.
Consider an embedding $\iota$ from a K3 surface $M$ to its 3-rd Hilbert scheme
of points $M^{[3]}$, mapping a point $x\in M$
to the subscheme with ideal sheaf $\calo_M / {\goth m}_x ^2$, where 
${\goth m}_x\subset \calo_M$ is the ideal of $x$.
It was shown in \cite{_Verbitsky:Hilbert_}.
that $\iota$ is not compatible with {\it any}
hyperk\"ahler structure on $M^{[3]}$. A number of similar
symplectic immersions can be easily constructed, but it was 
shown that for a generic K3 surface $M$, the Hilbert scheme
$M^{[n]}$ does not contain non-trivial hyperk\"ahler submanifolds.
This series of examples makes it natural to ask the following
question.

\hfill

\question
Let $\phi:\; X\hookrightarrow M$ be a holomorphic symplectic
immersion. Is there any way to establish whether $\phi$
is compatible with a hyperk\"ahler structure, using the holomorphic
symplectic geometry of $X$ and $M$?

\hfill

Let $M$ be a compact holomorphic symplectic 
K\"ahler manifold endowed with the unique hyperk\"{a}hler structure compatible
with the holomorphic symplectic structure and the K\"{a}hler class
(see \ref{_hyperka_symple_unique_Remark_}). 
In the present paper, we associate to any complex
subvariety $X\subset M$ a non-negative real number, called the
symplectic Wirtinger number $W(X)$ (\ref{_Wirtinger_numbe_Definition_}). 
We show that $W(X)\leq 1$, and that the equality is reached
if and only if $X$ is compatible with the hyperk\"ahler
structure on $M$ (\ref{_Wirtinger_Proposition_}). 
For a sequence of holomorphic symplectic immersions
\[ X_1\hookrightarrow X_2 \hookrightarrow X_3 ... \hookrightarrow M
\]
($\dim H^{2,0}(X_i) =1$, $\dim H^1(X_i)=0$ for all $i>1$, $\dim X_1>0$),
we prove the following inequality of Wirtinger numbers:
\[ 
W(X_1) \leq W(X_2) \leq W(X_3) \leq ...
\]
(\ref{_Wirti_increase_Theorem_})
In particular, if for some $k$ 
the submanifold $X_k$ is compatible with the hyperk\"ahler
structure, then for all $j>k$, the manifold $X_j$ is also
compatible with the hyperk\"ahler structure 
(\ref{_seque_Corollary_}).


\section{Hyperk\"ahler manifolds}


This subsection gives some of the basic and well known results 
and definitions from hyperk\"ahler geometry, most of which can be
found in \cite{_Besse:Einst_Manifo_} and in \cite{_Beauville_}.

\hfill

\definition \label{_hyperkahler_manifold_Definition_} 
(\cite{_Besse:Einst_Manifo_}) A {\bf hyperk\"ahler manifold} is a
Riemannian manifold $M$ endowed with three complex structures $I$, $J$
and $K$, such that the following holds.
 
\begin{description}
\item[(i)]  the metric on $M$ is K\"ahler with respect to these complex 
structures and
 
\item[(ii)] $I$, $J$ and $K$, considered as  endomorphisms
of the real tangent bundle, satisfy the relation 
$I\circ J=-J\circ I = K$.
\end{description}

\hfill

Clearly, a hyperk\"ahler manifold has a natural action of
the quaternion algebra ${\Bbb H}$ on its real tangent bundle $TM$. 
Therefore its complex dimension is even.
For each quaternion $L\in \Bbb H$, $L^2=-1$,
the corresponding automorphism of $TM$ is an almost complex
structure. It is easy to check that this almost 
complex structure is integrable (\cite{_Besse:Einst_Manifo_}).

\hfill

\definition \label{_indu_comple_str_Definition_} 
Let $M$ be a hyperk\"ahler manifold, and $L$ a quaternion satisfying
$L^2=-1$. The corresponding complex structure 
on $M$ is called {\bf an induced complex structure}. 
The $M$, considered as a K\"ahler manifold, is denoted by $(M, L)$. 
In this case, the hyperk\"ahler structure is called {\bf compatible
with the complex structure $L$}.

We say that the K\"ahler metric on $(M,L)$ is {\bf hyperk\"ahler}
if it is compatible with a hyperk\"ahler structure.

\hfill

Let $M$ be a hyperk\"ahler manifold; denote the
Riemannian form on $M$ by $(\cdot,\cdot)$.
Let the form $\omega_I := (\cdot, I(\cdot))$ be the usual K\"ahler
form  which is closed and parallel
(with respect to the Levi-Civita connection). Analogously defined 
forms $\omega_J$ and $\omega_K$ are
also closed and parallel. 
 
A simple linear algebraic
consideration (\cite{_Besse:Einst_Manifo_}) shows that the form
\begin{equation}\label{_holomo_symple_Equation_}
 \Omega:=\omega_J+\sqrt{-1}\omega_K 
\end{equation}
is of
type $(2,0)$ and, being closed, it is also holomorphic.
In addition, 
the form $\Omega$ is nowhere degenerate, as another linear 
algebraic argument shows.
It is called {\bf the canonical holomorphic symplectic form
of the manifold M}. Thus, for each hyperk\"ahler manifold $M$,
and an induced complex structure $L$, the underlying complex manifold
$(M,L)$ is holomorphically symplectic. 
Calabi's conjecture, proved by Yau in \cite{_Yau:Calabi-Yau_}, 
gives the converse statement.

\hfill

\theorem \label{_Calabi_Yau_Theorem_}
(\cite{_Beauville_},\cite{_Besse:Einst_Manifo_})
Let $X$ be a compact complex manifold
equipped with a holomorphic symplectic form, and let
$\omega$ be an arbitrary K\"ahler form on $X$. Then there exists a
unique hyperk\"ahler metric on $X$
with the K\"ahler form cohomologous to
 $\omega$. \endproof

\hfill

\remark\label{_hyperka_symple_unique_Remark_}
The metric does not determine uniquely the hyperk\"ahler structure:
there might by many hyperk\"ahler structures compatible with a givem
metric. However, from \ref{_holomo_symple_Equation_} it is clear
that the metrics together with the holomorphic symplectic form determine
the hyperk\"ahler structure uniquely.

\hfill

\definition \label{_simple_hyperkahler_mfolds_Definition_} 
(\cite{_Beauville_}) A connected simply connected 
compact hy\-per\-k\"ah\-ler manifold 
$M$ is called {\bf simple} if $M$ cannot be represented 
as a product of two hyperk\"ahler manifolds:
\[ 
   M\neq M_1\times M_2,\ \text{where} \ dim\; M_1>0 \ \ \text{and} 
   \ dim\; M_2>0
\] 
Bogomolov proved that every compact hyperk\"ahler manifold has a finite
covering which is a product of a compact torus
and several simple hyperk\"ahler manifolds. 
Bogomolov's theorem implies the following result
(\cite{_Beauville_}):

\hfill

\theorem\label{_simple_mani_crite_Theorem_}
Let $M$ be a compact hyperk\"ahler manifold.
Then the following conditions are equivalent.
\begin{description}
\item[(i)] $M$ is simple
\item[(ii)] $M$ satisfies $H^1(M, \R) =0$, $H^{2,0}(M) =\C$,
where $H^{2,0}(M)$ is the space of $(2,0)$-classes taken with
respect to some induced complex structure.
\end{description}

\endproof

\hfill

Calabi-Yau theorem can be stated in a more precise fashion, using
\ref{_simple_mani_crite_Theorem_}.

\hfill

\theorem \label{_symplectic_=>_hyperkahler_Proposition_}
(\cite{_Beauville_}, \cite{_Besse:Einst_Manifo_})
Let $M$, $\dim_\C M=n$, be a compact holomorphic
symplectic K\"ahler manifold with 
the holomorphic symplectic form $\Omega$, a K\"ahler class 
$[\omega]\in H^{1,1}(M)$ and a complex structure $I$. 
Assume that $\dim H^{2,0}(M) =1$, $H^1(M, \R) =0$, and
\begin{equation}\label{_inte_equa_Equation_}
\int_M \omega^n = \int_M (Re \Omega)^n.
\end{equation}
Then there exists a unique hyperk\"ahler 
structure $(I,J,K,(\cdot,\cdot))$
on $M$ such that the cohomology class of the symplectic form
$\omega_I=(\cdot,I\cdot)$ is equal to $[\omega]$ and the
canonical symplectic form $\omega_J+\1\omega_K$ is
equal to $\Omega$.

\hfill

{\bf Proof:} By Calabi-Yau (\ref{_Calabi_Yau_Theorem_}),
there exists a unique hyperk\"ahler metric on $M$ with the K\"ahler
class equal to $[\omega]$. Let $\c H_1=(I, J_1, K_1)$ be the corresponding 
hyperk\"ahler structure, and $\Omega_1:=\omega_{J_1}+\1 \omega_{K_1}$ 
its holomorphic symplectic form. Since $\dim H^{2,0}(M) =1$,
there exists a number $\lambda\in \C$ such that
$\Omega_1 =\lambda \Omega$.
A simple calculation shows that
\[ \int_M \omega^n =\int_M (Re \Omega_1)^n = \frac{1}{2^{n/2}}\binom{n}{n/2}
\int_M (\Omega_1\wedge \bar\Omega_1)^{\frac{n}{2}}
\]
(see \ref{omega_J,K_linear-alge_Lemma_})
and 
\[ \int_M (Re \Omega)^n = \frac{1}{2^{n/2}}\binom{n}{n/2}
\int_M (\Omega\wedge \bar\Omega)^{\frac{n}{2}}.
\]
{}From \eqref{_inte_equa_Equation_}, we obtain
\[ \int_M (\Omega\wedge \bar\Omega)^{\frac{n}{2}} = 
\int_M (\Omega_1\wedge \bar\Omega_1)^{\frac{n}{2}}.
\]
Therefore, $|\lambda| =1$. Write $\lambda$ in form $a + b \1$,
$a^2 + b^2=1$. Consider a hyperk\"ahler structure
$\c H$ with the same metrics, and quiaternionic action
given by the triple $(I, J, K)$: $J= a J_1 + b K$,
$K=-b J_1 +a K_1$. Clearly, the corresponding
form $\omega_J + \1 \Omega_K$ is equal to $\Omega$.
This proves \ref{_symplectic_=>_hyperkahler_Proposition_}.
\endproof

\hfill

\definition\label{_trianalytic_Definition_} 
(\cite{_Verbitsky:Symplectic_II_})
Let $X\subset M$ be a closed subset of a hyperk\"ahler
manifold $M$. Then $X$ is
called {\bf trianalytic} if $X$ is a complex analytic subset 
of $(M,L)$ for every induced complex structure $L$.
 
\hfill

Trianalytic subvarieties were a subject of a 
long study. Most importantly, consider a generic
induced complex structure $L$ on $M$. Then all closed complex
subvarieties of $(M, L)$ are trianalytic. 
Moreover,  a trianalytic subvariety can be canonically
desingularized (\cite{_Verbitsky:hypercomple_}),
and this desingularization is hyperk\"ahler.


\section{Wirtinger numbers and trianalytic subvarieties}
\label{_Wirti_defi_Section_}


Let $M$ be a compact complex manifold equipped with a 
hyperk\"ahler structure $\c H=(I, J, K, (\cdot, \cdot))$
and $X\subset M$ a closed complex subvariety
of even dimension. 
As usually, we denote by 
$\omega =\omega_I =(\cdot, I\cdot)$
the K\"ahler form of $M$, 
by 
\[ \Omega = \omega_J +\1 \omega_K = (\cdot, J\cdot)+\1 (\cdot, K\cdot)
\]
the holomorphic symplectic form, and
by $\bar \Omega$ the complex adjoint form $\omega_J -\1 \omega_K$.

Later on, we shall need the following lemma.

\hfill

\lemma \label{omega_J,K_linear-alge_Lemma_}
Let $M$ be a compact complex manifold equipped with a 
hyperk\"ahler structure and $X\subset M$ a closed complex subvariety
of even dimension $d$. 
Then 
\begin{equation}\label{_inte_omega_J,K_and_Omega_Lemma_}
   \binom{d}{\frac{1}{2}d} 
   \int_X \left(\frac{\Omega\wedge \bar \Omega}{2}\right)^{\frac{d}{2}}
   = \int_X \omega_J^{d} =
   \int_X \omega_K^{d}
\end{equation}

{\bf Proof:} The proof is pure linear algebra.
By definition, $\omega_J= \frac{1}{2}(\Omega+\bar \Omega)$,
and $\omega_K =-\frac{\1}{2}(\Omega-\bar \Omega)$,
Therefore,
\[ \int_X \omega_J^{d} = 
   \frac{1}{2^{d}}\sum_i \binom{d}{i}
\int_X \Omega^i\wedge\bar \Omega^{d-i}
\]
Since $X$ is complex analytic, $\int_X\eta=0$ unless 
$\eta$ is of type $(d, d)$. Therefore,
\[ \int_X \omega_J^{d} = 
   \frac{1}{2^{d}}\binom{d}{\frac{1}{2}d} 
\int_X \Omega^{\frac{d}{2}}\wedge\bar
   \Omega^{\frac{d}{2}}.
\]
The proof of the second equation is analogous.
\endproof

\hfill

Consider the natural $SU(2)$-action on $H^*(M)$.
Since the multiplication on cohomology of $M$ is $SU(2)$-invariant,
we have
\[ \int_M \omega_I^{\dim_\C M} = \int_M \omega_J^{\dim_\C M} =
   \int_M \omega_K^{\dim_\C M}.
\]
By \ref{omega_J,K_linear-alge_Lemma_},
\[ \int_M \omega_K^{\dim_\C M}=\int_M \omega_J^{\dim_\C M}
   =    \binom{\dim_\C M}{\frac{1}{2}\dim_\C M} 
   \int_M \left(\frac{\Omega\wedge \bar \Omega}{2}\right)^{\frac{\dim_\C M}{2}}.
\]
Therefore,
\[ \int_M \omega^{\dim_\C M} = 
   \binom{\dim_\C M}{\frac{1}{2}\dim_\C M} 
   \int_M \left(\frac{\Omega\wedge \bar \Omega}{2}\right)^{\frac{\dim_\C M}{2}}.
\]
By the same reasoning,
for any trianalytic subvariety $X\subset M$, we have
\[ \int_X \omega^{\dim_\C X} = 
      \binom{\dim_\C X}{\frac{1}{2}\dim_\C X} 
\int_X \left(\frac{\Omega\wedge \bar \Omega}{2}\right)^{\frac{\dim_\C X}{2}}.
\]
For an arbitrary complex subvariety $X\subset M$ of even dimension, 
consider the numbers
\[ \deg_\omega X:= \int_X \omega^{\dim_\C X} \]
and
\[ \deg_\Omega X:=    \binom{\dim_\C X}{\frac{1}{2}\dim_\C X} 
   \int_X \left(\frac{\Omega\wedge \bar \Omega}{4}\right)^{\frac{\dim_\C X}{2}}
\]
\definition
The number $\deg_\omega X$ is called {\bf the K\"ahler degree}
of $X$, and $\deg_\Omega X$ is called {\bf the symplectic degree}
of $X$.

\hfill

\remark\label{_Riemannian_Remark_}
Notice that 
\[ \deg_\omega X = ({\dim_\C X})! \Vol(X),
\]
where $\Vol(X)$ denotes the volume of $X$ taken with respect to the
Riemannian structure on $X\subset M$.

\hfill

The following fundamental inequality lays the groundwork
for all manipulations with trianalytic subvarieties.
Its K\"ahler version, the original Wirtinger's
inequality, is well known (see, e.g. \cite{_Stolzenberg_}).

\hfill

\proposition\label{_Wirtinger_Proposition_}
{\bf (Wirtinger inequality in holomorphic symplectic setting)}
Let $M$ be a compact complex manifold equipped with a 
hyperk\"ahler structure and $X\subset M$ a closed complex subvariety
of even dimension $d$. Then the following inequality holds
\begin{equation} \label{_Wirtinger_Equation_}
 \deg_\omega X \geq |\deg_\Omega X|.
\end{equation}
Moreover, \eqref{_Wirtinger_Equation_} is strict unless
$X$ is trianalytic.

\hfill

{\bf Proof:} 
This proof is essentially contained in \cite{_Verbitsky:Symplectic_II_}.
Let $\deg_{\omega_J} X$ denote the integral
\[ \int_X \omega_J^{\dim_\C X}.
\]
By \ref{omega_J,K_linear-alge_Lemma_},
\[ \deg_\Omega X = \deg_{\omega_J} X.
\]
Therefore, \eqref{_Wirtinger_Equation_}
can be written in the form
\[ |\deg_{\omega_J} X| \leq \deg_\omega X.\]
By \ref{_Riemannian_Remark_}, 
\[ \deg_\omega X = 2^d\Vol(X)
\]
By the classical Wirtinger inequality (\cite{_Stolzenberg_}, page 7),
\[ |\deg_{\omega_J} X|\leq  d!\Vol(X),
\]
and the equality is reached only if the subset
$X\subset M$ is complex analytic with respect
to $J$. This proves \eqref{_Wirtinger_Equation_}.
Finally, if \eqref{_Wirtinger_Equation_} is not strict,
then $X$ is complex analytic with respect to $J$ {\it and}
$I$. It is easy to show that such subvarieties are trianalytic.
\endproof

\hfill

\definition \label{_Wirtinger_numbe_Definition_}
In assumptions of \ref{_Wirtinger_Proposition_},
we define {\bf the Wir\-tinger's number} $W(X)$ of the subvariety
$X\subset M$ as 
\begin{equation}\label{_Wirte_defini_Equation_}
W(X):= \sqrt[d]{\frac{|\deg_\Omega(X)|}{\deg_\omega(X)}}.
\end{equation}
\ref{_Wirtinger_Proposition_} shows that
$W(X)\leq 1$, and
this inequality is strict $X$ unless is trianalytic.

\hfill

\remark \label{_Wirti_for_all_Remark_}
Notice that \eqref{_Wirte_defini_Equation_}
can be used to define the Wirtinger number for any
compact holomorphically symplectic K\"ahler manifold.


\section{Wirtinger numbers and symplectic immersions}


Let $M$ be a complex manifold. We say that $M$ is 
{\bf holomorphically symplectic} if $M$ is equipped
with nowhere degenerate holomorphic symplectic form $\Omega_M$.

\hfill

\definition
Let $\phi:\; X\hookrightarrow M$ be an immersion of holomorphic symplectic
manifolds. The map $\phi$ is called {\bf a holomorphic symplectic immersion}
if $\phi^* \Omega_M =\Omega_X$, where $\phi^*$ denotes the pullback
of differential forms.

\hfill

Wirtinger numbers, defined in Section 
\ref{_Wirti_defi_Section_},
give an interesting invariant
of holomorphic symplectic immersions. 

\hfill

\theorem \label{_Wirti_increase_Theorem_}
Let 
\[ 
   X_1\hookrightarrow X_2 \hookrightarrow ... \hookrightarrow X_n
\]
be a sequence of holomorphic symplectic immersions
of K\"ahler manifolds, $\dim X_1>0$. 
Assume that all manifolds $X_i$ are compact and K\"ahler, and the K\"ahler structures
are compatible with the immersions. Assume also that
for all $i>1$, the manifold $X_i$ is simple. Let $W(X_i)$ be the
Wirtinger numbers of $X_i$, in the sense of \ref{_Wirti_for_all_Remark_}.
Then
\[ 
  W(X_1) \leq W(X_2)\leq ...\leq W(X_n).
\]

\hfill

\corollary\label{_seque_Corollary_}
Let $M$ be a compact hyperk\"ahler manifold,
and $X\subset M$ a holomorphic symplectic subvariety which is simple.
Assume that $X$ contains a non-trivial subvariety which
is trianalytic in $M$. Then $X$ is trianalytic.

\hfill

{\bf Proof:}
Denote the trianalytic subvariety of $X$ by $Y$.
By \ref{_Wirtinger_Proposition_}, $W(Y)=W(M)$.
By \ref{_Wirti_increase_Theorem_}, 
$W(Y)\leq W(X)\leq W(M)$. Therefore,
$W(Y)= W(X)= W(M)$. Applying
\ref{_Wirtinger_Proposition_} again,
we obtain that $X$ is trianalytic.
\endproof

\hfill

{\bf Proof:} The proof of 
\ref{_Wirti_increase_Theorem_} uses the methods of hyperk\"ahler
geometry. Let $X_1\hookrightarrow X_2$
be a holomorphic symplectic immersion. It suffices
to show that $W(X_1) \leq W(X_2)$. From the definition
of Wirtinger numbers, it follows that 
\begin{equation}\label{_inte_Wirti_Equation_}
\deg_\Omega X_2 = \deg_{\omega'} X_2,
\end{equation}
where $\omega'$ is the K\"ahler form on $X_2$ given by
$\omega' = W(X_2) \omega$.
By \ref{_symplectic_=>_hyperkahler_Proposition_}, the
equality \eqref{_inte_Wirti_Equation_} implies that
there exists a hyperk\"ahler structure on $X_2$ such that
$[\omega']$ is its K\"ahler class and $\Omega$ its holomorphic
symplectic form. Using 
\ref{_Wirtinger_Proposition_},
we obtain that 
\begin{equation}\label{_Wirtinger_equ_for_symple_Equation_}
|\deg_\Omega X_1|\leq \deg_{\omega'} X_1. 
\end{equation} 
On the other hand, by the definition of $\omega'$, we have 
\begin{equation}\label{_Wirti_expre_Equation_}
   \deg_{\omega'} X_1 = \deg_{\omega} X_1\cdot W(X_2)^{\dim_\C X_1}.
\end{equation}
Dividing both sides of \eqref{_Wirtinger_equ_for_symple_Equation_}
by $\deg_{\omega} X_1$ and using \eqref{_Wirti_expre_Equation_},
we obtain
\[ W(X_2)^{\dim_\C X_1}\geq \frac{\deg_\Omega X_1 }{\deg_{\omega} X_1} 
   = W(X_1)^{\dim_\C X_1} 
\]
This proves \ref{_Wirti_increase_Theorem_}.
\endproof

\hfill

\hfill

{\bf Acknowledgements:}
I am grateful to D. Kaledin and A. Kuznetsov for insightful comments.

\hfill


\end{document}